\documentclass[graybox]{svmult}

\usepackage{makeidx}         
\usepackage{graphicx}        
\usepackage{multicol}        
\usepackage[bottom]{footmisc}

\usepackage{newtxtext}       %
\usepackage[varvw]{newtxmath}       

\makeindex             

\usepackage[colorlinks=true]{hyperref}
\newcommand{\gauss}{Gau{\ss}~}

\begin{document}

\title*{Groups, drift and harmonic measures}
\author{ Mark Pollicott and Polina Vytnova  }
\institute{M. Pollicott \at Department of Mathematics, Warwick University, Coventry,
   CV4 7AL, UK \\ \email{masdbl@warwick.ac.uk}
\and P. Vytnova \at Department of Mathematics, Warwick University, Coventry,
   CV4 7AL, UK \\ \email{P.Vytnova@warwick.ac.uk}}
%
%
\maketitle

\abstract*{
    In this short note we will describe an old problem and a new approach which
    casts light upon it.  The old problem is to understand the nature of
    harmonic measures for cocompact Fuchsian groups.  The new approach is to
    compute numerically the value of the drift and, in particular, 
    get new results on the dimension of the measure in some new examples.
}


\section{Introduction}
\label{sec:1}
An intriguing problem in modern geometric measure theory is the study of 
the harmonic measure on the unit circle which arises from a random walk on
a Fuchsian group. The existing approach combines several areas of pure
mathematics, such as  Ergodic Theory, Probability Theory, Hyperbolic Geometry, 
and rigorous Numerical Analysis.

We will first recall some background. Afterwards we introduce one of the quantifiable
characteristics of random walks called the \emph{drift} and explain how it is 
related to properties of the harmonic measure, in particular, its Hausdorff dimension.
Finally, we will draw a connection to a popular conjecture of Kaimanovich and Le Prince 
on the nature of the harmonic measure associated to a random walk on a Fuchsian group. 

Although there is a number of partial results in special cases the general
conjecture still remains open. We will offer a new perspective which covers both  some known examples and some new cases.
Finally, we will illustrate the question using the example of a
$(4,4,4)$-triangle group which can be traced back to the works of Gauss from 1805. 

\begin{quotation} 
    Wahrlich\footnote{``It is not knowledge, but the act of learning, not
    possession but the act of getting there, which grants the greatest
    enjoyment.''} es ist nicht das Wissen, sondern das Lernen, nicht das Besitzen
sondern das Erwerben, nicht das Da-Seyn, sondern das Hinkommen, was den
gr\"{o}ssten Genuss gew\"{a}hrt. 

\hfill Carl Friedrich \gauss to Wolfgang Bolyai~\cite{GtoB09}.
\end{quotation}

\section{Preliminaries}
In this section we collect together some background knowledge we need to properly formulate the problem.
\subsection{Hyperbolic Geometry}
We will treat Fuchsian groups as groups of isometries acting on the
hyperbolic plane~$\mathbb H$. For our considerations it will be convenient to consider
the so-called Poincar\'e disk model of~$\mathbb H$. 
\begin{figure}
    \sidecaption[t]
\includegraphics{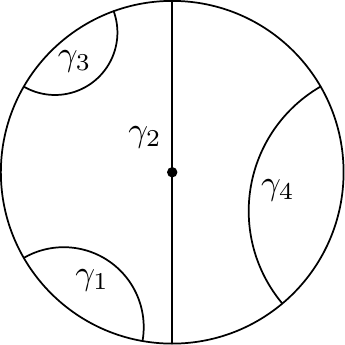}
\caption{The Poincar\'e disk where geodesics are either circular arcs which meet the
boundary circle orthogonally or diameters.}
\label{fig:1}
\end{figure}
\noindent This is a representation of the hyperbolic plane
as an open unit disk 
$\mathbb D = \{z = x+iy \colon |z| < 1\}$
equipped with the Poincar\'e metric
$$
ds^2 = 4 \frac{dx^2 + dy^2}{(1 - (x^2+y^2))^2}
$$
The geodesics are the extrema of the distance functional with respect to this
metric. They are precisely the Euclidean diameters and circular arcs which are orthogonal to
the boundary circle~$\partial \mathbb D$, as shown in Figure~\ref{fig:1}.

The orientation preserving isometries of~$\mathbb H$ in this model are 
linear fractional transformations of the form 
$$
g(z) = \frac{az + b}{\bar b z + \bar a}, \quad
\mbox{ where } a, b \in  \mathbb C \ \mbox{ and } \ |a|^2 - |b|^2 =1.
$$
As a group the orientation preserving isometries are isomorphic to the group 
consisting of~$2\times 2$ real matrices with determinant~$1$ up to multiplication 
by~$\pm \left(\begin{smallmatrix} 1 & 0 \\ 0 & 1 \end{smallmatrix}\right)$.
    Namely this is the group $\mathrm{PSL}(2, \mathbb R) = \mathrm{SL}(2, \mathbb R)/\{\pm I\}$.

\subsection{Geometric Group Theory}
In the present paper we want to consider finitely generated groups of
isometries of~$\mathbb H$. 
We call the group~$\Gamma$ {\it non-elementary} if it is not isomorphic to $\mathbb Z$. 
If $\Gamma$ is finitely generated, then the orbit $\Gamma 0 = \{g0 \colon g \in \Gamma\}$ 
of~$0 \in \mathbb D$ is a countable set of points in the unit disk. 
\begin{definition}
    We call a non-elementary group~$\Gamma$ a {\it Fuchsian group} if~$\Gamma 0$ is
    a discrete set with respect to the Poincar\'e metric~$ds^2$ introduced above. 
\end{definition}
In particular, if~$\Gamma$ is Fuchsian then all accumulation points with respect
to the Poincar\'e metric must lie on the unit circle $\partial \mathbb D =
\{z\in  \mathbb C \colon |z|=1\}$. 

Let us denote the set consisting of generators and their inverses
by~$\Gamma_\bullet=
\{g_1^{\pm1}, \cdots, g_d^{\pm 1}\}$.   
We can associate to~$\Gamma_\bullet$ the Cayley graph of~$\Gamma$.  
This is an infinite graph in which the vertices can be realised as the points 
of $\Gamma 0 = \{g0 \colon g \in \Gamma\}$ and two vertices~$g0$ and~$h0$ 
are connected by an edge
if and only if $gh^{-1} \in \Gamma_\bullet$. 

\subsection{Random walks and the drift}
We can now introduce our main tool. Given a set of generators and their
inverses~$\Gamma_\bullet$ we can consider a random walk on the Cayley graph 
where we allow a transition from a vertex~$g0$ to a neighbouring vertex~$h0$
with probability~$\frac{1}{2d}$. (Here we assume that $\#\Gamma_\bullet = 2d$,
as above.)

\begin{definition} 
    Given a specific set of generators~$\Gamma_\bullet$ we can associate the drift
    (or the rate of escape) defined by 
$$
\ell = \ell(\Gamma_\bullet) = \lim_{n\to +\infty} \frac{1}{(2d)^n} \sum_{g_{j_1}, \ldots, g_{j_n} \in \Gamma_\bullet}
\frac{d(g_{j_1} \cdots g_{j_n} 0,0)}{n}.
$$
\end{definition}
The limit always exists by a standard subadditivity argument and quantifies the rate at which typical points
$g_{j_1} \cdots g_{j_n} 0$ escape towards the boundary circle~$\partial \mathbb D$.

\section{Some examples}
Let us now turn to specific examples  of Fuchsian groups. For basic results on
hyperbolic polygons we refer the reader to an excellent book by
Beardon~\cite{beardon}.

\subsection{Regular octagon tilings}
\begin{figure}
    \centering
   \centerline{ \includegraphics{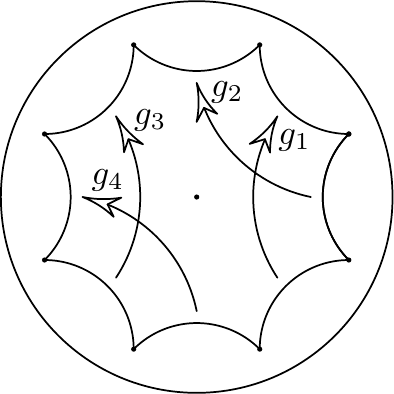}  \qquad \qquad
   \includegraphics{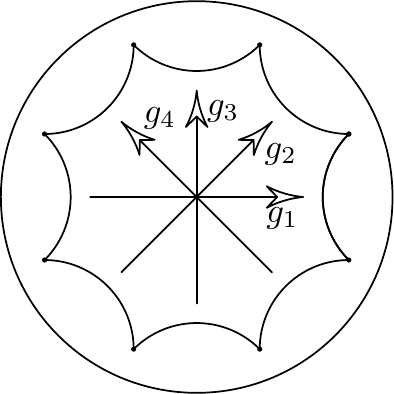} } 
\caption{ Left: Generators of the Gutzwiller group identify alternating sides of
the regular octagon; Right: Bolza group generated by isometries
identifying opposite sides of the regular octagon.} 
\label{fig:2}
\end{figure}

Let~$P \subset \mathbb H$ be a regular octagon with angles~$\frac{\pi}{4}$ and sides of equal
length which are geodesics as shown in Figure~\ref{fig:2}. We can consider groups of isometries
generated by four transformations which identify the sides of the regular octagon.
In this case the images of the octagon under the group action tile the
hyperbolic plane, so that $\Gamma P = \mathbb H$. 
There are four different identifications which yield a surface of
genus~$2$ as factor space $\mathbb H/\Gamma$~\cite{KN97}. This property
implies, in particular, that the group~$\Gamma$ generated by these identifications is
discrete. We will consider only two identifications which lead to well-known
surfaces: the Bolza surface and the Gutzwiller surface. 

\emph{The Gutzwiller group}~\cite{N95} $\Gamma_G = \langle g_1, g_2, g_3, g_4\rangle$ is generated 
by four isometries which identify the opposite sides of the regular pentagon~$P$. 
They satisfy the identity $g_1 g_2^{-1}g_1 g_2^{-1}g_3g_4^{-1} g_4 g_3^{-1} = I$
(see Figure~\ref{fig:2}, Left).

\emph{The Bolza group}~\cite{bolza} $\Gamma_B = \langle g_1, g_2, g_3, g_4\rangle$ is generated 
by four isometries which identify the opposite sides of the regular pentagon~$P$. 
They satisfy the identity $g_1 g_2^{-1}g_3 g_4^{-1}g_1^{-1}g_2 g_3^{-1}g_4 = I$
(see Figure~\ref{fig:2}, Right).

It turns out the drift doesn't depend on the identification chosen, in
particular, these two examples share the same value for the drift~$\ell$.

\begin{theorem}
    \label{thm:Bdrift}
With the choice of generators specified above, the drift~$\ell$ for the Bolza group~$\Gamma_B$
and for the Gutzwiller group~$\Gamma_G$ is the same and satisfies
$$
1.690771 < \ell < 1.691313.
$$
\end{theorem}
The method we use for estimating~$\ell$ involves looking at the action of~$\Gamma$ 
on~$\partial \mathbb D$ and computing the maximal Lyapunov exponent.
This is achieved by obtaining estimates on the spectral radius of transfer operators
acting on the space of $\alpha$-H\"older continuous functions 
$\mathcal L_t: C^\alpha(\partial \mathbb D) \to  C^\alpha(\partial \mathbb D)$
defined by $[\mathcal L_t f](z) = \frac{1}{2n}
\sum_{g \in \Gamma_\bullet}|g'(z)|^t f(gz)$
for~$t$ close to~$1$ and suitably small $\alpha > 0$. 
A more detailed exposition of the technical computer-assisted argument will
appear elsewhere. 

\subsection{Hyperbolic triangle groups}
Another class of interesting examples is perhaps the class of Coxeter groups
generated by reflections in the sides of a hyperbolic triangle, the so-called
triangle groups. It is easy to see that the group is discrete if and only if all angles of the triangle are
rational multipliers of~$\pi$. We will restrict our considerations to the case
when triangle has angles~$\frac{\pi}{k}$, $\frac{\pi}{l}$ and~$\frac{\pi}{m}$
where~$k,l$ and~$m$ are integers. 
Recall that sum of the angles of the hyperbolic triangle is strictly less
than~$\pi$, and therefore $k,l,m$ should satisfy the inequality $\frac{1}{k}+ \frac{1}{l}+\frac{1}{m} < 1$.

 \begin{figure}
     \centerline{\includegraphics{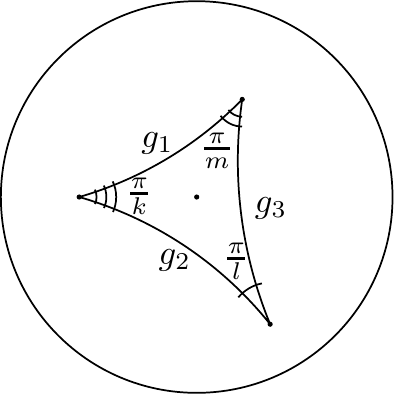} \qquad \qquad
     \includegraphics{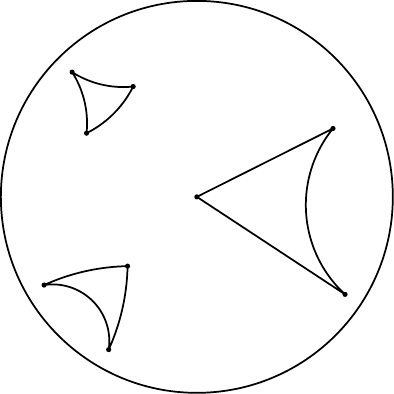}
}
\caption{Left: a triangle in~$\mathbb D$ with boundaries which are geodesics with
respect to the Poincar\'e metric and internal angles~$\frac{\pi}{k}$,
$\frac{\pi}{l}$ and~$\frac{\pi}{m}$, containing the centre of the disk in its
interior. Right: hyperbolic triangles with different
shape. }
\label{fig:3}
\end{figure}

\begin{definition}
The $(k,l,m)$-triangle group is a group generated by reflections in the 
three sides of a hyperbolic triangle with angles~$\frac{\pi}k$, $\frac{\pi}l$ and~$\frac{\pi}m$.
\end{definition}
It follows from the properties of reflections, that the generators~$g_1,g_2,g_3$
of the $(k,l,m)$-triangle group, shown in Figure~\ref{fig:3} satisfy the following relations
\begin{equation*}
g_1^2 = g_2^2 = g_3^2 = I \quad \mbox{ and } \quad (g_1g_2)^k = (g_2g_3)^l = (g_3g_1)^m = I.
\end{equation*}
These are the defining relations in the sense that any group with three generations
which satisfy these condition is the $(k,l,m)$-triangle group. Evidently, this
group is also cocompact. Furthermore, similarly to the case of the regular
octagon, the images of the original triangle with respect to the group form a
tessellation of the hyperbolic plane.

\begin{figure}
\centerline{
\includegraphics[width=53mm,height=61mm]{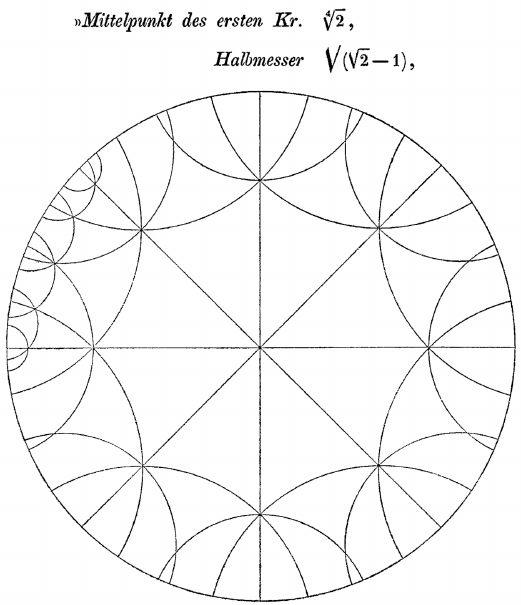} \qquad
\includegraphics[width=52mm]{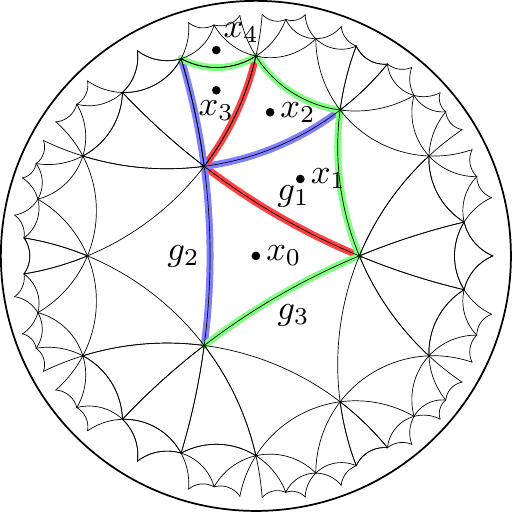}
}
\caption{Left: An original drawing by Gauss of a part of the $(4,4,4)$-triangle group
tessellation; Right: An orbit of a random walk corresponding to the
$(4,4,4)$-triangle group: $x_1 = g_1 x_0$, $x_2 = g_1 g_2 x0$, $x_3 = g_1 g_2
g_1 x_0$, and finally $x_4 = g_1 g_2 g_1 g_3 x_0$.  }
\label{fig:4}
\end{figure}

The study of the groups generated by reflections with respect to the sides of 
curvylinear triangles can be traced back to the works of Gauss. 
Bolyai, commenting on the Gauss' work, suggests that in a drawing
from ``Cereri Palladi Junoni Sacrum'' dated February 1805 Gauss introduced the
idea of reflection with respect to the circle. A copy of the drawing, 
taken from~\cite[p.~104]{Gauss} is shown in
Figure~\ref{fig:4} on the left. On the right we see a tessellation of the
hyperbolic plane generated by $(4,4,4)$-triangle group. The difference
between the two drawings is due to the choice of the location of the original
triangle. In the Gauss' drawing the centre of the disk is one of the vertices.
In our drawing, the centre of the disk is the barycentre of the triangle. 
Despite the appearance of a tessellation of the Poincar\'e disk in the
Gauss' drawing it was written 49 years before the birth of Poincar\'e! 
Lobachevsky laid the foundations of the hyperbolic geometry in 1823. 

Using the same machinery as in the case of the regular octagon we can estimate
the drift. In Table~\ref{tab:1} we list different examples of triangle
groups and give upper and lower bounds on the  associated drift. In a special 
case of the $(4,4,4)$-triangle group, we have the following result.
\begin{theorem}
    \label{thm:Tdrift}
The drift of the random walk on the $(4,4,4)$-triangle group with the choice of
generators specified above satisfies 
$$
0.1282273 < \ell < 0.1282264.
$$
\end{theorem}

\begin{table}[!t]
    \centering
\caption{Upper and lower bounds on the drift for some $(k,l,m)$-hyperbolic
triangle groups.}
\label{tab:1}
%
%
\begin{tabular}{p{2.5mm}cccp{5mm}p{28mm}p{5mm}p{28mm}}
\hline\noalign{\smallskip}
& k & l & m & & lower bound on $\ell$ & & upper bound on $\ell$ \\
\noalign{\smallskip}\svhline\noalign{\smallskip}
&  $ 3$&$7 $&$  2  $& & $    0.009936413804542  $& & $  0.009974294432083$ \\  
&  $ 3$&$8 $&$  2  $& & $    0.016242376981342  $& & $  0.016295700460901$ \\
&  $ 3$&$9 $&$  2  $& & $    0.020422904820936  $& & $  0.020508218335138$ \\
&  $ 4$&$5 $&$  2  $& & $    0.024263195172778  $& & $  0.024341830945392$ \\
&  $ 4$&$6 $&$  2  $& & $    0.037765501277040  $& & $  0.037870175386186$ \\
&  $ 4$&$8 $&$  2  $& & $    0.050724918174930  $& & $  0.050934249274956$ \\
&  $ 5$&$5 $&$  2  $& & $    0.046019792084900  $& & $  0.046155635941842$ \\
&  $ 5$&$6 $&$  2  $& & $    0.058159239428682  $& & $  0.058334985605960$ \\
&  $ 5$&$7 $&$  2  $& & $    0.065329026703739  $& & $  0.065563197936118$ \\
&  $ 6$&$6 $&$  2  $& & $    0.069559814745121  $& & $  0.069846131636394$ \\
&  $ 4$&$3 $&$  3  $& & $    0.046694831446660  $& & $  0.046816105401585$ \\
&  $ 5$&$3 $&$  3  $& & $    0.069435926662536  $& & $  0.069689191304812$ \\
&  $ 6$&$3 $&$  3  $& & $    0.081515978567027  $& & $  0.081925767935374$ \\
&  $ 7$&$3 $&$  3  $& & $    0.088431558608918  $& & $  0.089059709051931$ \\
&  $ 3$&$4 $&$  4  $& & $    0.088752444507380  $& & $  0.088919437571219$ \\
&  $ 3$&$6 $&$  6  $& & $    0.148515148139248  $& & $  0.149179933451390$ \\
&  $ 4$&$4 $&$  4  $& & $    0.128086862380309  $& & $  0.128344145942091$ \\
&  $ 5$&$5 $&$  5  $& & $    0.182618423778876  $& & $  0.183286144055414$ \\
&  $ 6$&$6 $&$  6  $& & $    0.209779208475952  $& & $  0.211031605163552$ \\
&  $ 7$&$7 $&$  7  $& & $    0.224864828238411  $& & $  0.228908301867331$ \\  
&  $ 8$&$8 $&$  8  $& & $    0.232248419011566  $& & $  0.238574707256068$ \\ 
&  $ 9$&$9 $&$  9  $& & $    0.236782098913020  $& & $  0.247054233672500$ \\
&  $10$&$10$&$  10 $& & $    0.240409132283172  $& & $  0.252180931190328$ \\
%
%
\noalign{\smallskip}\hline\noalign{\smallskip}
\end{tabular}
\end{table}


\section[Singularity]{Two problems}
After a brief discussion of the groups we are concerned with we continue by
introducing one of the central objects of the theory --- \emph{the harmonic
measure}.

\subsection{The harmonic measure on the unit circle}
As we have seen already, a typical orbit of the random walk associated to a
Fuchsian group converges to a point on the boundary circle 
with respect to the Euclidean and the hyperbolic metric. The distribution~$\nu$
of the limit points on the boundary~$\partial \mathbb D$ defines a
probability measure. More precisely, given a generating set~$\Gamma_\bullet$ of a
Fuchsian group we can define a family of 
probability measures on~$\mathbb D$ by 
$$
\nu_n = \left(\frac{1}{2d}\right)^n \sum_{g_{j_1} , \cdots, g_{j_n} \in \Gamma_\bullet}
\delta_{g_{j_1} \cdots g_{j_n}0}
$$
where $\delta_{g_{j_1} \cdots g_{j_n}0}$ is the Dirac measure supported at
$g_{j_1} \cdots g_{j_n}0$. The measures~$\nu_n$ converge in the weak star topology (on the closed unit
disk) to a probability measure~$\nu$ on~$\partial \mathbb D$.

\begin{definition}
    The measure~$\nu$ is called the \emph{harmonic measure} or the \emph{hitting measure}. 
\end{definition}

We can denote by~$\Lambda \subset \partial \mathbb D$ the support of this
measure (i.e., the smallest closed set of the full measure). It is known
that either $\Lambda = \partial \mathbb D$ or $\Lambda \subsetneq \partial
\mathbb D$ is a Cantor set. 

\subsection{Singularity of the harmonic measure}
The following natural question was posed by Kaimanovich and Le Prince~\cite{kp}: 

\smallskip
\noindent{\bf Question 1} Can we characterise Fuchsian groups for which the associate
harmonic measure is absolutely continuous with respect to the Lebesgue measure?

\smallskip


Of course, if the support of the harmonic measure is a Cantor set then the
measure is  singular with respect to Lebesgue measure. 
Therefore, we will only consider the case that $\Lambda = \partial \mathbb D$. 
In the special case when one of the generators in~$\Gamma_\bullet$ is
parabolic it was shown by Gadre, Maher, and Tiozzo that the harmonic measure is
always singular~\cite{gadre}. 

Furthermore, there are examples of non-discrete groups due to Bourgain 
for which~$\nu$ is absolutely continuous~\cite{bourgain} (see also~\cite{bps}).
In the more general setting when the weights in the random walk differ
the measure~$\nu$ may be singular~\cite{kp}.

  In the setting of the surface groups, this question has been intensively studied~\cite{kosenko},~\cite{kt}. 
 One set of examples is Fuchsian groups generated by isometries identifying 
 the sides of hyperbolic polygons. However, many of the examples with
 more than four sides have harmonic measures that are singular
 (see~\cite[Theorem~1]{kosenko}).

\subsection{Dimension of the harmonic measure}
An important quantitative characteristic of the measure is its Hausdorff
dimension. 

\begin{definition}
    The Hausdorff dimension of the measure is the infimum of Hausdorff
    dimensions of sets of the full measure:  
$$
\dim_H(\nu) = \inf\{\dim_H(X) \mid X \subset \partial \mathbb D \mbox{ Borel
and } \nu(X) = 1 \}.  
$$
\end{definition}
There is a useful result due to Tanaka which relates the question of absolute
continuity of the harmonic measure to the numerical value of its Hausdorff
dimension~\cite{tanaka}. 

\begin{proposition}[Tanaka] 
    \label{prop:tanaka}
    The harmonic measure~$\nu$ is absolutely continuous if and only if~$\dim_H(\nu)=1$.
\end{proposition}

This leads to the following stronger version of Question~$1$.

\smallskip

\noindent{\bf Question 2} Assuming the harmonic measure is not absolutely continuous with respect to
Lebesgue measure, can we estimate its Hausdorff dimension? 

\smallskip

We now return to our examples.

  It follows from the result of Kosenko~\cite[Theorem 1.2]{kosenko} (see
  also~\cite{kt}) that the harmonic measure~$\nu_B$ associated to \emph{the
  Bolza group}~$\Gamma_B$ is singular. 
  We can improve on this result. 

  \begin{theorem} The dimension of the harmonic measure~$\nu_B$ for the Bolza group
      satisfies 
  $$
     \dim_H(\nu_B) \leq  0.86116.
  $$
  The dimension of the harmonic measure~$\nu_G$ for the Gutzwiller group~$\Gamma_G$       satisfies
  $$
        \dim_H(\nu_G) \leq 0.86317.
  $$
  In particular, the  harmonic measure~$\nu_G$ is also singular. 
  \label{thm:dimnu}
  \end{theorem}  

  In order to explain the proof, we need one extra ingredient. 

\subsection{Relation to the Avez entropy}
We introduce another numerical characteristic of the random walk which is
commonly used to estimate the dimension of a measure.
\begin{definition} 
 We can associate to a harmonic measure~$\nu$ the {\it Avez random walk entropy} defined by 
 $$
 h_A(\nu) = \lim_{n \to +\infty} \frac{1}{n} H(\nu^{* n})
 $$
 where~$H(\cdot)$ is the usual Shannon entropy function and~$\nu^{* n}$ denotes the
 $n$-fold convolution~\cite{A74}.
\end{definition}

The limit always exists by subadditivity.  
The dimension, entropy and drift are related  using the following
identity~\cite{bhm},~\cite{GMM15},~\cite{tanaka}.

\begin{proposition}
    \label{prop:gouezel}
    For the harmonic measure~$\nu$ we have that $\dim_H(\nu) = h_A(\nu)/\ell(\nu)$.
\end{proposition}
Now we are ready to prove Theorem~\ref{thm:dimnu}. \\
{\it Proof of Theorem~\ref{thm:dimnu}.}
Combining Proposition~\ref{prop:tanaka} and Proposition~\ref{prop:gouezel} we
see that in order to establish that the harmonic  measure is  singular it is
sufficient to show that $h(\nu) \geq \ell(\nu)$. In
particular, we need to establish an upper bound on the entropy and a lower bound
on the drift. For the Gutzwiller group an estimate on the entropy is given in a
beautiful paper~\cite{GMM15}. Even the most basic bound they give in Example~2.3
$h(\nu_B)\le 1.46$ in combination with the estimate on the
drift~$\ell(\nu_B)$ from Theorem~\ref{thm:Bdrift} allows us to deduce that the
measure in singular.    
In the case of the Bolza group, we can use the estimate on the
drift~$\ell(\nu_G)$ from Theorem~\ref{thm:Bdrift} and the upper bound 
 on the entropy of~$h(\nu_G)\le \frac{3}{4}\log 7$ coming from the free group on four generators.

 \subsection{Final remarks}
Should one wish to apply the same approach to show that the harmonic measure for the triangle 
groups is singular it will be necessary to obtain an effective upper bound on the Avez entropy.   
Unfortunately, the naive bound of $\frac13{\log 2} \approx 0.231049\ldots $
corresponding to the Avez entropy of the random walk on the free product 
$\mathbb Z_2 \times \mathbb Z_2 \times \mathbb Z_2$ isn't quite low enough to
show that the measure is singular for most of triangle groups listed in Table~\ref{tab:1}. 
Nevertheless in the case of $(8,8,8)$-,$(9,9,9)$-, and $(10,10,10)$-triangle
groups we may conclude that the measure is singular. It is reasonable to suggest
that the drift for $(k,k,k)$-triangle group is monotone increasing as $k \to \infty$. This
would imply that the harmonic measure is singular for $k\ge 8$. 

\begin{acknowledgement}
 The first author is partly
supported by ERC-Advanced Grant 833802-Resonances and EPSRC grant
EP/T001674/1 the second author is partly supported by
EPSRC grant EP/T001674/1.
\end{acknowledgement}


\end{document}